\documentclass[10pt,a4paper]{article}
\usepackage{amsmath,amssymb,amsthm}
\usepackage[pdftex,colorlinks = true, citecolor = green, urlcolor = black]{hyperref}
\begin{document}

\newtheorem{theorem}{Thm}

\title{Maximum Probability  and Relative Entropy Maximization.
Bayesian Maximum Probability  and Empirical Likelihood}
\author{M. Grend\'ar\footnote{Dept. of Mathematics, Bel University,
Tajovskeho 40, 974 01 Banska Bystrica, Slovakia. E-mail:
marian.grendar@savba.sk. Inst. of Measurement Science, SAS,
Bratislava. Inst. of Mathematics and CS, SAS, Banska Bystrica. Date:
Apr 9, 2008. To appear in Proc. of Intnl. Workshop on Applied
Probability (IWAP) 2008, Compi\`egne, France, July 7-10, 2008.}}
\date{}
\maketitle

\begin{abstract}
Works, briefly surveyed here, are concerned with two basic methods:
Maximum Probability   and Bayesian Maximum Probability; as well as
with their asymptotic instances: Relative Entropy Maximization   and
Maximum Non-paramet\-ric Likelihood. Parametric and empirical
extensions of the latter methods -- Empirical Maximum Maximum
Entropy and Empirical Likelihood  -- are also mentioned. The methods
are viewed as tools for solving certain ill-posed inverse problems,
called $\Pi$-problem, $\Phi$-problem, respectively. Within the two
classes of problems, probabilistic justification and interpretation
of the
respective methods are  discussed.\\

\end{abstract}

\section{$\Phi$-problem, MAP,  MNPL}

The $\Phi$-problem can be loosely stated as follows: there is a
prior distribution over a non-parametric set $\Phi$ of data-sampling
distributions and a sample from unknown data-sampling distribution.
The objective is to select a data-sampling distribution from the set
$\Phi$, called model.

More formally: Let $\mathcal P$ be the set of all probability mass
functions\footnote {For the sake of simplicity the presentation is
restricted to the discrete case.  The continuous case is treated in
\cite{GJa}.} (pmf's) with finite support $\mathcal X$. The set
$\mathcal{P}$ is endowed with the usual topology. Let
$\Phi\subseteq\mathcal{P}$. Let $X_1^n \triangleq X_1, X_2, \dots,
X_n$ be i.i.d. sample from pmf $r\in\mathcal P$. The 'true' sampling
distribution $r$ need not be in $\Phi$; in other words: the model
$\Phi$ might be misspecified. A strictly positive prior $\pi(\cdot)$
is put over $\Phi$. The objective in the $\Phi$-problem is to select
a sampling distribution $q$ from $\Phi$, when the information
summarized by $\{\mathcal X, X_1^n, \pi(\cdot), \Phi\}$ and nothing
else is available.

Bayesian Maximum Probability method selects the Maximum A-Posteriori
Probable (MAP) data-sampling distribution(s)  $\hat q_\textrm{MAP}
\triangleq \arg\sup_{q\in\Phi} \pi_n(q\,|\,X_1^n)$;  there the
posterior probability $\pi_n(q|X_1^n) \propto {e^{-l_n(q)}\pi(q)}$,
and $l_n(q)$ is used to denote $-\sum_{i=1}^n \log q(x_i)$; $\log$
is meant with the base $e$. Hence the standard abbreviation, MAP,
for the method.

The Bayesian Sanov Theorem (BST), through its corollary -- the
Bayesian Law of Large Numbers (BLLN) -- provides a strong case for
MAP as the correct method for solving the $\Phi$-problem. The
theorems are Bayesian counterparts of the well-known Large
Deviations (LD) theorems for empirical measures: the Sanov Theorem
and the Conditional Law of Large Numbers (cf. \cite{C} and Sect. 2).
In order to state the theorems it is necessary to introduce the
$L$-divergence $L(q||p)$ of $q\in\mathcal P$ with respect to
$p\in\mathcal P$: $L(q||p) \triangleq - \sum_\mathcal X p \log q$.
The {$L$-projection} $\hat q$ of $p$ on $Q\subseteq\mathcal P$ is
$\hat q \triangleq \arg \inf_{q \in Q} L(q||p)$. The value of
$L$-divergence at an $L$-projection of $p$ on $Q$ is denoted by
$L(Q||p)$.

\begin{theorem} {\rm (BST)} Let $X_1^n$ be i.i.d. $r$. Let $Q\subset\Phi\subseteq\mathcal{P}$;
$L(Q\,||\,r) < \infty$. Then for $n\rightarrow\infty$, $
\frac{1}{n}\log \pi_n(q\in Q|X_1^n) = - \{L(Q\,||\,r) -
L(\Phi\,||\,r)\}, {a.s.\ } r^\infty. $
\end{theorem}

The posterior probability $\pi_n(Q|X_1^n)$ decays exponentially fast
(a.s. $r^\infty$) with the decay rate $L(Q\,||\,r) -
L(\Phi\,||\,r)$. For a proof see \cite{GL}. To the best of our
knowledge Ben-Tal, Brown and Smith \cite{BBS} were the first to use
an LD reasoning in the Bayesian nonparametric setting. Ganesh and
O'Connell \cite{GO} proved BST for the well-specified special case;
i.e., $r\in\Phi$, by means of formal LD.

\begin{theorem} {\rm(BLLN)} Let $\Phi\subseteq\mathcal P$ be a convex,
closed set. Let $B(\hat q, \epsilon)$, be a closed $\epsilon$-ball
defined by the total variation metric, centered at the
$L$-projection $\hat q$ of $r$ on $\Phi$. Then, $
\lim_{n\rightarrow\infty}\pi_n(q\in B(\hat q,\epsilon)\,|\,X_1^n) =
1, {a.s.\ } r^\infty.$
\end{theorem}

The BLLN is an extension of Freedman's Bayesian nonparametric
consistency theorem \cite{F} to the case of misspecified model. It
shows that the posterior probability concentrates (a.s. $r^\infty$)
on the $L$-projection of the 'true' sampling distribution $r$ on
$\Phi$. For a book-length treatment of Bayesian non-parametric
consistency see \cite{GR}.

MAP satisfies the BLLN. To see this, note that by the Strong Law of
Large Numbers (SLLN), conditions for supremum of  the posterior
probability asymptotically turn into conditions for supremum of the
negative of $L$-divergence. This also permits to view the
$L$-projections as asymptotic instances of MAP distributions $\hat
q_\textrm{MAP}$.

There is also another method which satisfies the BLLN: Maximum
Non-parametric Likelihood (MNPL). This can be shown by the above
mentioned recourse to the SLLN. MNPL selects $\hat q_\textrm{MNPL}
\triangleq \arg\sup_{q\in Q} - l_n(q)$.

These two (up to trivial transformations) are the only methods for
solving the $\Phi$-problem, which comply with the BLLN; hence they
are consistent in the well-specified as well as in the misspecified
case. Selecting a sampling distribution by some other conceivable
method would, in general, asymptotically select sampling
distribution which is {\em a posteriori} zero-probable. In this
sense, selection of, say, the posterior mean, or selection of
$\arg\sup_{q\in\Phi} -\sum_{\mathcal X} q\log\frac{q}{r}$, are ruled
out.

The $\Phi$-problem becomes more interesting when turned into a
parametric setting. To this end, let $X$ be a random variable with
pmf $r(x; \theta)$ parame\-triz\-ed by $\theta\in\Theta\subseteq
\mathbb R^K$. Assume that a researcher is not willing to specify
parametric family $q(X; \theta)$ of data-sampling distributions, but
is only willing to specify some of its underlying features. These
features, i.e., the model $\Phi$, can be characterized by Estimating
Equations (EE): $\Phi\triangleq\bigcup_\Theta \Phi(\theta)$, where
$\Phi(\theta) \triangleq \{q(x; \theta): \sum_{\mathcal X}
q(x;\theta)u_j(x;\theta) = 0, 1\le j\le J\}$,
$\theta\in\Theta\subseteq\mathbb R^K$. In the EE theory parlance,
$u(\cdot)$ are the estimating functions, number of which is in
general different than the number $K$ of parameters $\theta$. The
'true' data sampling distribution $r(x; \theta)$ need not belong to
$\Phi$. A Bayesian puts positive prior $\pi$ over $\Phi$, which in
turn induces prior $\pi(\theta)$ over $\Theta$; cf. \cite{FR}. By
the BLLN, the posterior $\pi_n(\cdot|X_1^n)$ concentrates on a weak
neighborhood of the $L$-projection $\hat q$ of $r(x;\theta)$ on
$\Phi$:
\begin{equation*}
\hat q(x;\hat\theta) =
\arg\inf_{\theta\in\Theta}\inf_{q(x;\theta)\in\Phi(\theta)}
L(q(x;\theta)\,||\,r(x;\theta)).
%-\int r(x)\log q(x;\theta).
\end{equation*}
This thus provides a probabilistic justification for using
$\hat\theta$ as an estimator of $\theta$. Thanks to the convex
duality, the estimator $\hat\theta$ can be obtained also as
$\hat\theta =
\arg\sup_{\theta\in\Theta}\inf_{\lambda(\theta)\in\mathbb R^J}
-\sum_{i=1}^m r(x_i)\log(1 - \sum_j \lambda_j(\theta)
u_j(x_i;\theta))$. Since $r$ is in practice not known, following
\cite{KS}, one can estimate the convex dual objective function by
$-\sum_{l=1}^n \log(1 - \sum_j \lambda_j(\theta) u_j(x_l;\theta))$.
The resulting estimator is just the Empirical Likelihood (EL)
estimator (cf. \cite{QL}, \cite{O}, \cite{MJM}). It can be easily
seen that EL satisfies the BLLN. The same is true for the Bayesian
MAP estimator $ \hat q_{\textrm{MAP}}(x;\hat\theta_{\textrm{MAP}}) =
\arg\sup_{\theta\in\Theta}\sup_{q(x;\theta)\in\Phi(\theta)}
\pi_n(q(x;\theta)\,|\,X_1^n).$
For further results and discussion see \cite{GJa}, \cite{GJb}.

\section{$\Pi$-problem, MaxProb, REM}

Unlike the $\Phi$ problem, the $\Pi$ problem is not a statistical
problem. In the $\Pi$ problem, the sampling distribution $q$ is
known, and there is a set $\Pi\subseteq\mathcal P$, into which an
unavailable empirical pmf, drawn from $q$, is assumed to belong. The
objective is to select an empirical pmf (also known as type, cf.
\cite{C}) from the set $\Pi$. Thus, the $\Phi$ and $\Pi$ problems
are opposite to each other.

More formally: let $\mathcal X$ be a set of $m$ elements. Type
$\nu^n \triangleq [n_1, n_2, \dots, n_m]/n$, where $n_i$ is the
number of occurrences of $i$-th element of $\mathcal X$ (i.e.,
outcome), $i = 1,2,\dots, m$, in a sample of size $n$, drawn from
sampling distribution $q$. The objective in the $\Pi$-problem is to
select a type(s) $\nu^n$ from $\Pi$, when the information summarized
by $\{\mathcal X, q, n, \Pi\}$ and nothing else is available.

Maximum Probability (MaxProb) method (cf. \cite{B}, \cite{Vi},
\cite{GG}) selects the type $\hat\nu^n = \arg\sup_{\nu^n\in\Pi}
\pi(\nu^n; q)$  which can be generated by the sampling distribution
$q$, with the highest probability. If the sampling is i.i.d., then
$\pi(\nu^n; q) = n\,!\prod_{i=1}^m\frac{q_i^{n_i}}{n_i\,!}$. Niven
\cite{N1} expanded MaxProb  into non-i.i.d. and combinatorial
settings; see also \cite{N2}, \cite{Vi}, \cite{GN}.

The Sanov Theorem (ST) (cf. \cite{S}, \cite{CST}), through its
corollary -- the Conditional Law of Large Numbers (CLLN) (cf.
\cite{Va}, \cite{CC}, \cite{CST}) -- provides a probabilistic
justification for application of MaxProb in the i.i.d. instance of
the $\Pi$-problem. The ST identifies the exponential decay rate
function as the $I$-divergence $I(p\,||\,q) \triangleq \sum
p\log\frac{p}{q}$, $p,q\in\mathcal P$. The {$I$-projection} $\hat p$
of $q$ on $\Pi\subseteq\mathcal P$ is $\hat p \triangleq \arg
\inf_{p \in \Pi} I(p\,||\,q)$. The value of the $I$-divergence at an
$I$-projection of $q$ on $\Pi$ is denoted by $I(\Pi||q)$.

\begin{theorem} {\rm (ST)}
 Let $\Pi$ be an open set; $I(\Pi\,||\,q)<\infty$. Then, for
 $n\rightarrow\infty$,
  $\frac{1}{n} \log \pi(\nu^n \in\Pi; q) = -
  I(\Pi\,||\,q).$
\end{theorem}

The rate of the exponential convergence of the probability
$\pi(\nu^n \in \Pi; q)$ towards zero is determined by the
information divergence at (any of) the $I$-projection(s) of $q$ on
$\Pi$.

\begin{theorem} {\rm (CLLN)}
Let $\Pi$ be a convex, closed set that does not contain $q$. Let
$B(\hat p, \epsilon)$ be a closed $\epsilon$-ball defined by the
total variation metric that is centered at the $I$-projection $\hat
p$ of $q$ on $\Pi$. Then, $\lim_{n \rightarrow \infty} \pi(\nu^n \in
B(\hat p, \epsilon)\,|\,\nu^n \in \Pi; q) = 1.$
\end{theorem}
Given that a type from $\Pi$ was observed, it is asymptotically
zero-probable that the type was different than the $I$-projection of
the sampling distribution $q$ on $\Pi$.

It is straightforward to see that MaxProb satisfies CLLN. Indeed,
set of MaxProb types converges to set of $I$-projections, as
$n\rightarrow\infty$; cf. \cite{GAI}, \cite{GG}. Relative Entropy
Maximization method (REM/MaxEnt) which maximizes, with respect to
$p$, the negative of $I$-divergence (a.k.a., relative entropy) thus
can be viewed as asymptotic form of MaxProb method.

Still, it is possible to solve $\Pi$-problem by selecting the
type(s) with the highest value of relative entropy; in other words,
to view REM as a self-standing method for solving $\Pi$-problem,
rather than as an asymptotic instance of MaxProb. Obviously, REM
satisfies CLLN.

MaxProb and REM/MaxEnt are the only two methods which satisfy CLLN.
Selection of the mean type, which was under the name ExpOc proposed
in \cite{GG}, or selection of, say the type with the highest value
of Tsallis entropy, would in general, violate CLLN.

The $\Pi$-problem originated in Statistical Physics, where $\Pi$ is
formed by mean energy constraint; see \cite{E}. In \cite{GGA}
feasible set of types formed by interval observations was
considered.

Estimating Equations can be used to expand the $\Pi$ problem into
parametric setting. This time, the EE define a feasible set $\Pi$
into which an unobserved parametrized type $\nu^n(\theta)$ is
supposed to belong: $\Pi \triangleq \bigcup_\Theta\Pi(\theta)$,
where $\Pi(\theta) \triangleq \{p(x;\theta): \sum_{\mathcal X}
p(x;\theta)u_j(x;\theta) = 0, 1\le j \le J\}$,
$\theta\in\Theta\subseteq \mathbb R^K$. The true data-sampling
distribution $r(x; \theta)$  need not belong to $\Pi$. The
parametric $\Pi$-problem is framed by the information $\{\mathcal X,
r, n, \Pi(\theta), \Theta\}$, and the objective is now to select
parametric type $\nu^n(\theta)$ from $\Pi$.
CLLN implies  (cf. \cite{KSE}) that the parametric $\Pi$-problem
should be (for $n\rightarrow\infty$) solved by selecting
$$
\hat p(x;\hat\theta) =
\arg\inf_{\theta\in\Theta}\inf_{p(x;\theta)\in\Pi(\theta)}
I(p(x;\theta)\,||\,r(x;\theta)).
$$
Thanks to the convex duality, the estimator $\hat\theta$ can
equivalently be obtained as $\hat\theta =
\arg\sup_{\theta\in\Theta}\inf_{\lambda(\theta)\in\mathbb R^J}
\log\sum_{i=1}^m r(x_i; \theta)\exp(-\sum_{j=1}^J \lambda_j(\theta)
u_j(x_i;\theta))$. The estimator is known as {Maximum Maximum
Entropy} (MaxMaxEnt) estimator.

The parametric $\Pi$-problem can be made more realistic, by assuming
that a sample of size $N$ is available to a modeler. Kitamura and
Stutzer \cite{KS} suggested to use the sample to estimate the convex
dual objective function by its sample analogue $\log\sum_{l=1}^N
\exp(-\sum_{j=1}^J \lambda_j u_j(x_l;\theta))$. The resulting method
is known as Empirical Maximum Maximum Entropy (EMME) method, or
Maximum Entropy Empirical Likelihood (cf. \cite{ISJ}, \cite{KS},
\cite{MJM}, \cite{JM}).

\section{Acknowledgements}

Valuable discussions with George Judge  and Robert Niven, and a
feedback from Val\'erie Girardin are gratefully acknowledged.
Supported by VEGA 1/3016/06 and APVV RPEU-0008-06 grants.

\end{document}